\documentclass[12pt]{article}
\usepackage{bbm}
\usepackage{latexsym}
\usepackage{mathrsfs}
\usepackage{graphicx}
\usepackage{subfigure}
\usepackage{amssymb}
\usepackage{amsmath}
\usepackage{amsthm}
\usepackage{color}
\usepackage{cite}
\usepackage{float}
\usepackage{indentfirst}
\usepackage{anysize}\marginsize{45mm}{45mm}{40mm}{50mm}
\usepackage{caption}
\usepackage{float}
\usepackage{multirow}

\usepackage{enumerate}

\usepackage{latexsym, bm}
\newtheoremstyle{lemma}{\topsep}{\topsep}%
     {}
     {}
     {\bfseries}
     {}
     {0.1em}
     {\thmname{#1}\thmnumber{ #2}\thmnote{ #3}}
\theoremstyle{lemma}  

\newtheorem{theorem}{Theorem}     
\newtheorem{lemma}[theorem]{Lemma}
\newtheorem{corollary}[theorem]{Corollary}
\newtheorem{conjecture}[theorem]{Conjecture}

\newtheorem{definition}{Definition}

\numberwithin{equation}{section}

\title{ On the conjecture of vertex-transitivity of Dcell\thanks{This research was partially supported by the fundamental research funds for the central universities (No. 2672018ZYGX2018J069)}}

\author{ Huazhong L\"{u}\\
{\small $^{1}$School of Mathematical Sciences,} \\
{\small University of Electronic Science and Technology of China,}\\
{\small Chengdu, Sichuan 610054, P.R. China}\\
{\small E-mail: lvhz08@lzu.edu.cn}\\}

\date{}
\begin{document}

\maketitle
\begin{abstract}

Gu et al. in [Inform. Process. Lett. 134 (2018) 52--56] conjectured that the data center network $D_{k,n}$ is vertex-transitive for all $k\geq0$ and $n\geq2$. In this paper, we show that $D_{k,n}$ is vertex-transitive for $k\leq1$ and $n\geq2$, and it is not vertex-transitive for all $k\geq2$ and $n\geq2$.

\vskip 0.1 in

\noindent \textbf{Key words:} Data center network; Vertex-transitivity
\end{abstract}

\section{Introduction}

Data centers are crucial to the business of companies such as Amazon, Google and Microsoft. Data centers with tens of thousands servers were built to offer increasingly popular on-line applications such as web search, email, cloud storage, on-line gaming, etc. Guo et al. \cite{Guo} proposed the data center networks, briefly say DCell, $D_{k,n}$ for parallel computing systems, which has numerous favourable features for data center networking. In DCell, a large number of servers are connected by high-speed links and switches, providing much higher network capacity compared the tree-based.

In many situations, as it often simplifies the computation and routing algorithms, parallel interconnect of high symmetry is highly desirable. To deal with the symmetry of a graph, the aim is to obtain as much information as possible about the symmetric property of the graph. Gu et al. \cite{Gu} showed that $D_{k,n}$ $(k\geq0,n\geq2)$ is not edge transitive except the cases of $k=0$ and $k=1,n=2$. In the same paper, by several base cases of $D_{k,n}$, they conjectured the following:

\begin{conjecture}{\bf.}
$D_{k,n}$ is vertex-transitive for all $k\geq0$ and $n\geq2$.
\end{conjecture}

\vskip 0.05 in

We solve this conjecture in this paper.

The rest of this paper is organized as follows. In Section 2, some notations and the definitions of Dcell are presented. The main results of this paper are shown in Section 3.

\vskip 0.05 in

\section{Preliminaries}

Let $G=(V(G),E(G))$ be a graph, where $V(G)$ is vertex-set of $G$ and $E(G)$ is edge-set of $G$. The number of vertices of $G$ is denoted by $|G|$. A path $P=\langle x_0,x_1,\cdots, x_k\rangle$ in $G$ is a sequence of distinct vertices so that there is an edge joining consecutive vertices. If a path $C=\langle x_0,x_1,\cdots, x_k\rangle$ is such that $k\geq3$, $x_0=x_k$, then $C$ is said to be a {\em cycle}, and the length of $C$ is the number of edges contained in $C$. For other standard graph notations not defined here please refer to \cite{Bondy}.

In what follows, we shall present the definition of the DCell.
\vskip 0.0 in

\begin{definition}{\bf .}\label{def}\cite{Wang}
A $k$ level DCell for each $k$ and some global constant $n$, denoted by $D_{k,n}$, is recursively defined as follows. Let $D_{0,n}$ be the complete graph $K_n$ and let $t_{k,n}$ be the number of vertices in $D_{k,n}$. For $k\geq1$, $D_{k,n}$ is constructed from $t_{k-1,n}+1$ disjoint copies of $D_{k-1,n}$, where $D_{k-1,n}^i$ denotes the $i$th copy. Each pair of $D_{k-1,n}^a$ and $D_{k-1,n}^b$ ($a<b$) is joined by a unique $k$ level edge below.

A vertex of $D_{k-1,n}^i$ is labeled by $(i,a_{k-1},\cdots,a_0)$, where $k\geq1$ and $a_0\in\{0,1,\cdots,$ $n-1\}$. The suffix $(a_{j},a_{j-1},\cdots,a_0)$, of a vertex $v$, has the unique $uid_j$, given by $uid_j(v)=a_0+\sum_{l=1}^{j}(a_lt_{l-1,n})$. The vertex $uid_{k-1}$ $b-1$ of $D_{k-1,n}^a$ is connected to $uid_{k-1}$ $a$ of $D_{k-1,n}^b$.
\end{definition}

\begin{figure}
\centering
\includegraphics[height=80mm]{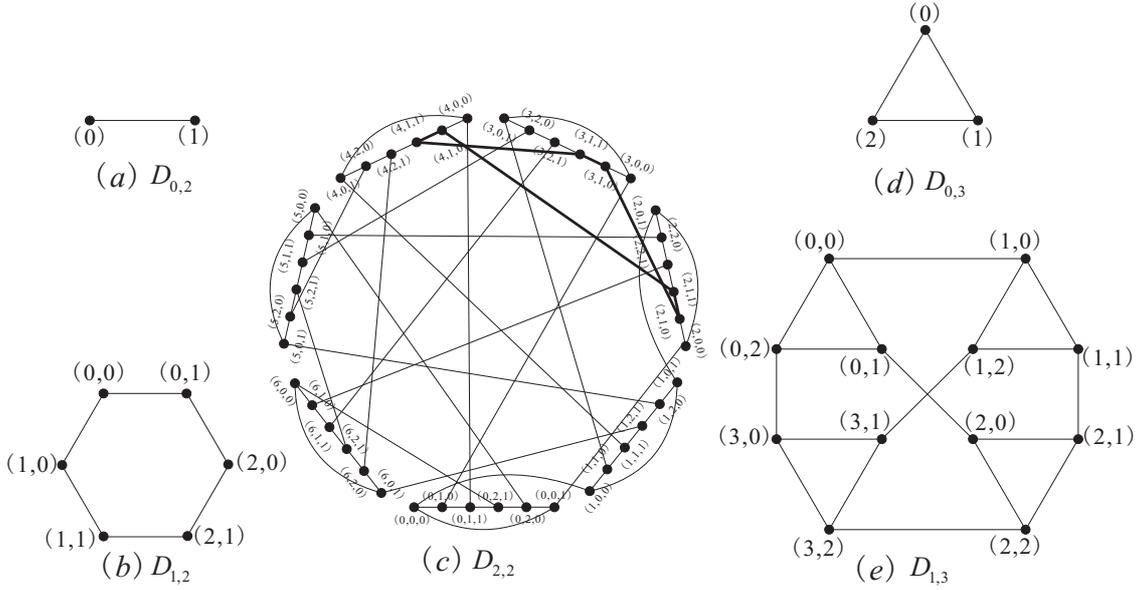}
\caption{Some small DCells.} \label{DCell}
\end{figure}

By definition above, it is obvious that $D_{0,2}$ is an edge, $D_{0,3}$ is a triangle and $D_{1,2}$ is a 6-cycle. Several $D_{k,n}$ with small parameters $k$ and $n$ are illustrated in Fig. \ref{DCell}. The edge between $D_{k-1,n}^a$ and $D_{k-1,n}^b$ is said to be a {\em level} $k$ edge. For convenience, let $E_k$ denote the set of all level $k$ edges of $D_{k,n}$.

\section{Main results}

Since $D_{0,n}$ is isomorphic to the complete graph $K_n$, it is vertex-transitive for all $n\geq2$. In what follows, we shall consider the vertex-transitivity of $D_{1,n}$ for all $n\geq2$.

Toward this end, we prove the following theorem which is clearly stronger than the above and so implies the result. Let $H$ be the graph containing $n+1$ disjoint copies of the complete graph $K_n$, $n\geq2$, for convenience, denoted by $K_n^i$, $i\in\{0,1,\cdots,n\}$. There is exact one edge between any two copies of $K_n^i$, which implies that each vertex in $K_n^i$ has exact one neighbor not in $K_n^i$.

\begin{theorem}\label{v-t-H}{\bf.} $H$ is vertex-transitive for all $n\geq2$.
\end{theorem}
\noindent{\bf Proof.} Our aim is to show that for any two vertices $u,v\in V(H)$, there exists an automorphism of $H$ that maps $u$ to $v$. By the definition $H$, we know that $H$ consists of $n+1$ copies of $K_n$ and each pairs of distinct $K_n$s is connected by an edge. For convenience, the vertices of $H$ are labeled as $(i,j)$, $0\leq i\leq n$, $0\leq j\leq n-1$, where $i$ means vertices in $K_n^i$. Suppose without loss of generality that $u=(i_1,j_1)$ and $v=(i_2,j_2)$. We consider the following two cases.

\noindent{\bf Case 1.} $i_1\neq i_2$. Then $u$ and $v$ lie in $K_{n}^{i_1}$ and $K_{n}^{i_2}$, respectively.

\noindent{\bf Case 1.1.} $uv\in E(H)$. We define a map $f:V(K_{n}^{i_1})\rightarrow V(K_{n}^{i_2})$ as follows, (1) $f(u)=v$; (2) if $x\in V(K_{n}^{i_1})$ and $y\in V(K_{n}^{i_2})$ have neighbors in the same $K_{n}^{j}$ ($j\neq i_1,i_2$), respectively, let $f(x)=y$. Thus, we can obtain a permutation $\varphi$ that swaps vertices of $K_{n}^{i_1}$ with that of $K_{n}^{i_2}$ under $f$, and fixes other vertices in $H$. Clearly, $\varphi$ preserves adjacency of $H$, indicating that $H$ is an automorphism with $\varphi(u)=v$.

\noindent{\bf Case 1.2.} $uv\not\in E(H)$. Let $u'$ (resp. $v'$) be a neighbor of $u$ (resp. $v$) not in $K_{n}^{i_1}$ (resp. $K_{n}^{i_2}$). By definition of $H$, there are three cases of relative positions of $u'$ and $v'$.

\noindent{\bf Case 1.2.1.} $u'\in V(K_{n}^{i_2})$ and $v'\in V(K_{n}^{i_3})$ ($i_3\neq i_1,i_2$). There exists an edge $w'w$ between $K_{n}^{i_1}$ and $K_{n}^{i_3}$, where $w'\in V(K_{n}^{i_1})$ and $w\in V(K_{n}^{i_3})$. Clearly, $w\neq v'$ and $w'\neq u$. We define an automorphism $\varphi$ such that, (1) $u\rightarrow v \rightarrow w\rightarrow u$ and $u'\rightarrow v'\rightarrow w'\rightarrow u'$; (2) if $x_1\in V(K_{n}^{i_1})$, $x_2\in V(K_{n}^{i_2})$ and $x_3\in V(K_{n}^{i_3})$ have neighbors $y_1,y_2$ and $y_3$ in the same $K_{n}^{j}$ ($j\neq i_1,i_2,i_3$), respectively, then $x_1\rightarrow x_2 \rightarrow x_3\rightarrow x_1$ and $y_1\rightarrow y_2\rightarrow y_3\rightarrow y_1$; (3) $\varphi$ fixes other vertices in $H$. Obviously, $\varphi$ preserves adjacency of $H$ such that $\varphi(u)=v$.

\noindent{\bf Case 1.2.2.} $u',v'\in V(K_{n}^{i_3})$ ($i_3\neq i_1,i_2$). Then there exists an edge $ww'$ between $K_{n}^{i_1}$ and $K_{n}^{i_2}$, where $w\in V(K_{n}^{i_1})$ and $w'\in V(K_{n}^{i_2})$. Clearly, $w\neq u$ and $w'\neq v$. There exists an automorphism $\varphi$ such that, (1) $\varphi(u)=v$, $\varphi(w)=w'$ and $\varphi(u')=v'$; (2) if $x_1\in V(K_{n}^{i_1})$ and $y_1\in V(K_{n}^{i_2})$ have neighbors in the same $K_{n}^{j}$ ($j\neq i_1,i_2,i_3$), respectively, then $\varphi(x_1)=y_1$; (3) $\varphi$ fixes all other vertices. Obviously, $\varphi$ preserves adjacency of $H$ such that $\varphi(u)=v$.

\noindent{\bf Case 1.2.3.} $u'\in V(K_{n}^{i_4})$ and $v'\in V(K_{n}^{i_3})$. ($i_3,i_4\neq i_1,i_2$). Then there exists an edge $ww'$ from $K_{n}^{i_1}$ to $K_{n}^{i_2}$. Similarly, there exist an edge $w_1'w_1$ from $K_{n}^{i_4}$ to $K_{n}^{i_3}$, an edge $u_1u_1'$ from $K_{n}^{i_1}$ to $K_{n}^{i_3}$, and an edge $v_1v_1'$ from $K_{n}^{i_2}$ to $K_{n}^{i_4}$. There exists an automorphism $\varphi$ such that, (1) $\varphi(u)=v$, $\varphi(w)=w'$, $\varphi(u')=v'$, $\varphi(w_1')=w_1$, $\varphi(u_1)=v_1$ and $\varphi(u_1')=v_1'$; (2) if $x_1\in V(K_{n}^{i_1})$, $x_2\in V(K_{n}^{i_2})$, $x_3\in V(K_{n}^{i_3})$ and $x_4\in V(K_{n}^{i_4})$ have neighbors $y_1$, $y_2$, $y_3$ and $y_4$ in the same $K_{n}^{j}$ ($j\neq i_1,i_2,i_3,i_4$), respectively, then $\varphi(x_1)=x_2$, $\varphi(x_3)=x_4$, $\varphi(y_1)=y_2$ and $\varphi(y_3)=y_4$. Moreover, let $\varphi$ fix other vertices in $H$. Obviously, $\varphi$ preserves adjacency of $H$ such that $\varphi(u)=v$.

\noindent{\bf Case 2.} $i_1=i_2$. Then $u$ and $v$ lie in $K_{n}^{i_1}$. We first choose a permutation $P$ on vertices of $K_{n}^{i_1}$ such that $P$ swaps $u$ with $v$ and fixes other vertices in $K_{n}^{i_1}$. Let $u_1=(i_3,j_3)$ and $v_1=(i_4,j_4)$ be neighbors of $u$ and $v$ not in $K_{n}^{i_1}$, respectively. Clearly, $i_3\neq i_4$. We define a map $f:V(K_{n}^{i_3})\rightarrow V(K_{n}^{i_4})$ as follows, (1) $f(u_1)=v_1$; (2) if $x\in V(K_{n}^{i_3})$ and $y\in V(K_{n}^{i_4})$ have neighbors in the same $K_{n}^{j}$ ($j\neq i_1,i_3,i_4$), respectively, let $f(x)=y$; (3) if $xy$ is an edge from $K_{n}^{i_3}$ to $K_{n}^{i_4}$, let $f(x)=y$. Obviously, $f$ is a bijection. So, we then swap vertices of $K_{n}^{i_3}$ with that of $K_{n}^{i_4}$ under $f$. Thus, we can obtain a permutation $\psi$ that swaps $u$ with $v$ in $K_{n}^{i_1}$, swaps vertices of $K_{n}^{i_3}$ with that of $K_{n}^{i_4}$ according to $f$, and fixes other vertices in $H$. Clearly, $\psi$ preserves adjacency of $H$, indicating that $H$ is an automorphism with $\psi(u)=v$. \qed

\vskip 0.1 in

It is obvious that $D_{1,n}$ is a special case of $H$, so the following corollary is straightforward.

\begin{corollary}\label{v-t-01-n}{\bf.} $D_{1,n}$ is vertex-transitive for all $n\geq2$.
\end{corollary}

Gu et al. (Conclusion in \cite{Gu}) pointed that $D_{2,2}$ is vertex-transitive by MAGMA software. However, this is not true. By detailed computation of MAGMA software, $D_{2,2}$ is not vertex-transitive. In fact, (0,2,0) is contained in exact one 6-cycle of $D_{2,2}$, namely $D_{1,2}^0$, while (3,1,1) is contained in at least two 6-cycles, namely $D_{1,2}^3$ and $\langle (3,1,1),(3,1,0),(2,1,0),(2,1,1),(4,1,0),(4,1,1),(3,1,1)\rangle$ (see heavy lines in Fig. \ref{DCell} $(c)$). Let $u=(0,\cdots,0,2,1)$ be a vertex in $D_{k,2}$, based on the fact above, we have the following lemma.

\begin{lemma}{\bf.}\label{one-6-cycle} $u$ is contained in exact one 6-cycle of $D_{k,2}$ for all $k\geq2$.
\end{lemma}

\noindent{\bf Proof.} We proceed by induction on $k$. It is known that $u$ is contained in exact one 6-cycle of $D_{2,2}$, thus, the induction basis holds. We assume that $u$ is contained in exact one 6-cycle of $D_{k-1,2}$ for $k\geq 3$. Next we consider $D_{k,2}$. Since $u\in V(D_{k-1,2}^0)$, by the induction hypothesis, $u$ is contained in exact one 6-cycle of $D_{k-1,2}^0$. The following vertices $u_i$, $0\leq i\leq k$, are level $i$ neighbors of $u$ in $D_{k,2}$.

$u_0=(0,\cdots,0,2,0)$,

$u_1=(0,\cdots,0,1,1)$,

$u_2=(0,\cdots,0,6,0,0)$,

$u_3=(0,\cdots,6,0,0,0)$,

$\cdots\cdots$

$u_{k-1}=(0,6,\cdots,0,0,0)$,

$u_k=(6,0,\cdots,0,0,0)$.

Observe that if contract each $D_{k-1,2}^j$ ($0\leq j\leq |D_{k-1,2}|$) as a single vertex, then the resulting graph is isomorphic to $K_{|D_{k-1,2}|+1}$. To form another 6-cycle $C$ of $D_{k,2}$ containing $u$, by structure of $D_{k,2}$, we have $|E(C)\cap E_k|=3$. Thus, $uu_k\in E(C)$. Noting $u_k\in D_{k-1,2}^6$, then $C$ contains exact one edge in each of $D_{k-1,2}^0$, $D_{k-1,2}^6$. So exact one of $uu_0, uu_1,\cdots,uu_{k-1}$ is contained in $C$. Noting also that there exists exact one edge between any two distinct $D_{k-1,2}^j$s, then each of $u_i$ (except $u_k$) has exact one neighbor $u_{i,k}$ not in $D_{k-1,2}^0$. Our aim is to verify all possible $C$ containing exact one of $uu_0, uu_1,\cdots,uu_{k-1}$. For brevity, we present $u_{i,k}$ as follows.

$u_{0,k}=(5,0,\cdots,0,0,0)$,

$u_{1,k}=(4,0,\cdots,0,0,0)$,

$u_{2,k}=(37,0,\cdots,0,0,0)$,

$u_{3,k}=(253,0,\cdots,0,0,0)$,

$\cdots\cdots$

$u_{k-1,k}=(6\times|D_{k-2,2}|+1,0,\cdots,0,0,0)$.

The edges $e_{6,l}$ between $D_{k-1,2}^6$ and $D_{k-1,2}^l$ are as follows, where $l\in\{4,5,\cdots,6\times|D_{k-2,2}|+1\}$.

$e_{6,4}=(6,0,\cdots,0,2,0)(4,0,\cdots,0,2,1)$,

$e_{6,5}=(6,0,\cdots,0,2,1)(5,0,\cdots,0,2,1)$,

$e_{6,37}=(6,0,\cdots,0,0,6,0,0)(37,0,\cdots,0,1,0,0)$,

$e_{6,253}=(6,0,\cdots,0,6,0,0,0)(253,0,\cdots,0,1,0,0)$,

$\cdots\cdots$

$e_{6,6\times|D_{k-1,2}|+1}=(6,6,0,\cdots,0,0,0)(6\times|D_{k-2,2}|+1,0,\cdots,0,1,0,0)$.

By checking, it is not hard to see that if one of $e_{6,l}$ is contained in $C$, then $C$ is not a 6-cycle. Thus, the lemma holds. \qed

\begin{lemma}{\bf.}\label{not-v-t-k-2} $D_{k,2}$ is not vertex-transitive for all $k\geq2$.
\end{lemma}
\noindent{\bf Proof.} By Lemma \ref{one-6-cycle}, it can be known that $u=(0,\cdots,0,0,2,1)$ is contained in exact one 6-cycle of $D_{k,2}$, while $v=(0,\cdots,0,3,1,1)$ is contained in at least two 6-cycles. Thus, the statement follows immediately. \qed

\begin{lemma}{\bf.}\label{not-v-t-k-3-n} $D_{k,n}$ is not vertex-transitive for all $k\geq2$ and $n\geq3$.
\end{lemma}
\noindent{\bf Proof.} For convenience, we denote the number of 6-cycles containing a vertex $v$ in $D_{k,n}$ by $c_n^k(v)$. Let $u=(0,\cdots,0,0,0)$ and $v=(0,\cdots,0,1,2)$ be two vertices in $D_{k,n}$. To prove this lemma, we shall show that $c_n^k(u)\neq c_n^1(u)$ and $c_n^k(v)=c_n^1(v)$ for all $k\geq2$, which clearly implies that $D_{k,n}$ is not vertex-transitive. For a given $n\geq3$, we proceed the lemma by induction on $k$. We first consider $k=2$. So, at this time, $u=(0,0,0)$ and $v=(0,1,2)$ are two vertices in $D_{1,n}^0$ of $D_{2,n}$. Since $D_{1,n}^0\cong D_{1,n}$, noting $D_{1,n}$ is vertex-transitive, we have $c_n^1(u)=c_n^1(v)$. Clearly, $C=\langle u,(0,0,1),(2,0,0),(2,0,1),(1,0,1),(1,0,0),u,\rangle$ is a 6-cycle in $D_{2,n}$ such that $V(C)\not\subseteq V(D_{1,n}^0)$. Thus, $c_n^k(u)>c_n^1(u)$ for all $k\geq2$. It remains to show that $c_n^k(v)=c_n^1(v)$ for all $k\geq2$. We need to show that there exists no 6-cycle $C'$ of $D_{2,n}$ containing $v$ such that $V(C')\not\subseteq V(D_{1,n}^0)$. Analogous to the proof of Lemma \ref{one-6-cycle}, it is not hard to see that $c_n^2(v)=c_n^1(v)$ holds. Thus, the induction basis holds. So we assume that $c_n^{k-1}(v)=c_n^1(v)$ for $k\geq3$. Next we consider $D_{k,n}$.

%
%
%
%
%
%
%
%
%
%
%
%
%
%
%
%
%
%
%
%
%
%
%
%

Let $v_0^i$, $0\leq i\leq n-1$, $i\neq2$, and $v_j$, $1\leq j\leq k$, be $0$ and $j$ level neighbors of $v$ in $D_{k,n}$, respectively. For convenience, we list them below.

$v_0^0=(0,\cdots,0,1,0)$,

$v_0^1=(0,\cdots,0,1,1)$,

$v_0^3=(0,\cdots,0,1,3)$,

$\cdots\cdots$

$v_0^{n-1}=(0,\cdots,0,1,n-1)$,

$v_1=(0,\cdots,0,3,1)$,

$v_2=(0,\cdots,n+3,0,0)$,

$\cdots\cdots$

$v_{k-1}=(0,n+3,0,\cdots,0,0,0)$,

$v_{k}=(n+3,0,\cdots,0,0,0)$.

Similarly, to form a 6-cycle $C'$ (not a subgraph of $D_{k-1,n}^0$) of $D_{k,n}$ containing $v$, it is known that $|E(C')\cap E_k|=3$. Thus, $vv_k\in E(C')$. Noting that $v_k\in D_{k-1,n}^{n+3}$, then $C'$ contains exact one edge of each of $D_{k-1,n}^0$ and $D_{k-1,n}^{n+3}$. So exact one of $vv_0^i$ and $vv_j$ is contained in $C'$. Noting also that there exists exact one edge between any two distinct $D_{k-1,n}$s in $D_{k,n}$, then each of $v_0^{i}$ and $v_j$ has exact one neighbor $v_{0,k}^i$ and $v_{j,k}$ not in $D_{k-1,n}^0$, respectively. In what follows, we shall verify if there possibly exists such a 6-cycle $C'$. For brevity, we present $u_{0,k}^i$ and $v_{j,k}$ as follows.

$v_{0,k}^0=(n+1,0,\cdots,0)$,

$v_{0,k}^1=(n+2,0,\cdots,0)$,

$v_{0,k}^3=(n+4,0,\cdots,0)$,

$\cdots\cdots$

$v_{0,k}^{n-1}=(2n,0,\cdots,0)$,

$v_{1,k}=(3n+2,0,\cdots,0)$,

$v_{2,k}=(n(n+1)(n+3)+1,0,\cdots,0)$,

$\cdots\cdots$

$v_{k-1,k}=((n+3)|D_{k-2,n}|+1,0,\cdots,0)$.

The edges $e_{n+3,l}$ between $D_{k-1,n}^{n+3}$ and $D_{1,n}^{l}$ are as follows, where $l\in\{n+1,n+2,n+4,\cdots,2n\}\cup\{3n+2,n(n+1)(n+3)+1,\cdots,(n+3)|D_{k-2,n}|+1\}$.

$e_{n+3,n+1}=(n+3,0,\cdots,0,1,1)(n+1,0,\cdots,0,1,2)$,

$e_{n+3,n+2}=(n+3,0,\cdots,0,1,2)(n+2,0,\cdots,0,1,2)$,

$e_{n+3,n+4}=(n+3,0,\cdots,0,1,3)(n+4,0,\cdots,0,1,3)$,

$\cdots\cdots$

$e_{n+3,2n}=(n+3,0,\cdots,0,1,n-1)(2n,0,\cdots,0,1,3)$,

$e_{n+3,3n+2}=(n+3,0,\cdots,0,3,1)(3n+2,0,\cdots,0,1,3)$.

$e_{n+3,n(n+1)(n+3)+1}=(n+3,0,\cdots,0,n+3,0,0)(n(n+1)(n+3)+1,0,\cdots,0,1,3)$,

$\cdots\cdots$

$e_{n+3,(n+3)|D_{k-2,n}|+1}=(n+3,n+3,0,\cdots,0,0,0)((n+3)|D_{k-2,n}|+1,0,\cdots,0,1,3)$.

By checking, it is not hard to see that if one of $e_{n+3,l}$ is contained in $C'$, then $C'$ is not a 6-cycle, which implies that $c_n^{k}(v)=c_n^1(v)$. (Noting that when $n=3$, the last two coordinates of vertices with ($*$,1,3) will degenerated to ($*$,2,2), where ``$*$'' denotes the first coordinate. The statement is also true.) Thus, the lemma holds. \qed

Combining Lemmas \ref{not-v-t-k-2} and \ref{not-v-t-k-3-n}, the following theorem is straightforward.

\begin{theorem}{\bf.}\label{not-v-t-k-n} $D_{k,n}$ is not vertex-transitive for all $k\geq2$ and $n\geq2$.
\end{theorem}

%

%

\end{document}